# On Helical Surfaces with a Constant Ratio of Principal Curvatures

YANG LIU[1], OLIMJONI PIRAHMAD[2], HUI WANG[2], DOMINIK L. MICHELS[1], HELMUT POTTMANN[2],

[1]Computational Sciences Group,

[2]Computational Design and Fabrication Research Group,

KAUST Visual Computing Center

We determine all helical surfaces in three-dimensional Euclidean space which possess a constant ratio $a := \kappa_1/\kappa_2$ of principal curvatures (CRPC surfaces), thus providing the first explicit CRPC surfaces beyond the known rotational ones. A key ingredient in the successful determination of these surfaces is the proper choice of generating profiles. We employ the contours for parallel projection orthogonal to the helical axis. This has the advantage that the CRPC property can be nicely expressed with the help of the involution of conjugate surface tangents. The arising ordinary differential equation has an explicit parametric solution, which forms the basis for a further study and classification of the possible shapes and the singularities arising for $a > 0$.

## 1 INTRODUCTION

Surfaces which possess a relation $F(\kappa_1, \kappa_2) = 0$ between their principal curvatures $\kappa_1$ and $\kappa_2$ are usually named after the German mathematician Julius Weingarten [17] who – back in the 19th century – studied them in connection with a characterization of those surfaces that are isometric to rotational surfaces. The latter are exactly the focal surfaces of Weingarten surfaces. There has been a huge interest in very special cases, such as surfaces with constant Gaussian or mean curvature, but apart from those, there are only very few explicitly known Weingarten surfaces.

Apart from the important theoretical aspects, Weingarten surfaces are also interesting for the realization of freeform facades in architecture: Due to the functional relation between principal curvatures, a Weingarten surface possesses only a one-parameter family of curvature elements. This facilitates cost effective paneling solutions due to the possibility of working with a number of molds that is significantly smaller than the number of panels [8, 9].

In this paper, we contribute to surfaces $S$ in three-dimensional Euclidean space which possess a constant ratio of principal curvatures $\kappa_1/\kappa_2 =: a$. We exclude surfaces with one vanishing principal curvature (developable surfaces) and call the others *CRPC surfaces*. The original motivation for our interest in these surfaces also arises from architectural geometry, since negatively curved CRPC surfaces appear in so-called asymptotic gridshells with a constant node angle [12, 13]. CRPC surfaces are also interesting shapes for the so-called Caravel meshes of Tellier et al. [15]. There, the principal curvature ratio appears in the aspect ratio of faces in an asymptotic analysis of Caravel quad meshes and Caravel hexagonal meshes with edge offsets, so that CRPC surfaces give rise to visually well balanced architectural structures.

CRPC surfaces are a natural generalization of minimal surfaces ($a = -1$). However, in big contrast to minimal surfaces, very little is known about them. Explicit parameterizations are only available for rotational CRPC surfaces [1, 3, 5–7, 16]. Lopez and Pampano [5] recently presented a classification of all surfaces with a linear relation $\kappa_1 = a\kappa_2 + b$ between principal curvatures, including a study of the special case $b = 0$ of rotational CRPC surfaces. Their work also contains

Author's address: Yang Liu[1], Olimjoni Pirahmad[2], Hui Wang[2], Dominik L. Michels[1], Helmut Pottmann[2], {yang.liu.4,pirahmad.olimjoni,hui.wang.1, dominik.michels,helmut.pottmann}@kaust.edu.sa,
[1]Computational Sciences Group,
[2]Computational Design and Fabrication Research Group,
KAUST Visual Computing Center.

a variational characterization of the profiles of these surfaces. Rotational CRPC surfaces with $K < 0$ have also been characterized via isogonal asymptotic parameterizations $f(u,v)$ where $\|f_u\| = \|f_v\|$ [10, 11, 14]. Moreover, it has been shown that Weingarten surfaces to a linear relation of the form $a\kappa_1 + b\kappa_2 + c = 0$ are rotational if they are foliated by a family of circles [4].

Jimenez et al. [2] derived CRPC surfaces via a Christoffel-type transformation of certain spherical nets, with a focus on discrete models. Effective methods for the computation of discrete CRPC surfaces [16] provided some insight on the shape variety of CRPC surfaces. As these are based on numerical optimization, one cannot derive precise mathematical conclusions, but conjectures as basis for further studies.

To add new explicit representations of CRPC surfaces beyond the rotational ones, one will look into other special surface classes. Unfortunately, some classes are excluded quickly: A *ruled CRPC surface* has the rulings as one family of asymptotic curves and the other family of asymptotic curves needed to intersect the rulings under a constant angle. The set of second asymptotic directions $A(t)$ along a ruling $R$ is a regulus on a ruled quadric. This contradicts a constant angle between $R$ and $A(t)$ except for a right angle (leading the helicoid as ruled minimal surface). One can also apply a result by Beltrami and Dini [3], which states that ruled Weingarten surfaces are helical or rotational, and check that there are no CRPC surfaces among them except helicoids. Likewise, channel surfaces (envelopes of families of spheres) can only be Weingarten surfaces if they are rotational surfaces or helical pipe surfaces, so that the only CRPC channel surfaces are just the rotational ones.

In this paper, we explicitly determine all *helical CRPC surfaces*. Clearly, this amounts to solving a 2nd order ordinary differential equation (ODE). However, this ODE may turn out so complicated that an explicit solution appears to be very difficult. We show how to choose proper generating curves of the helical surface so that the ODE is sufficiently simple to be solved explicitly. In order to come up with that equation, we do not compute the principal curvatures, but work only with the involutory projective map of conjugate surface tangents.

An idea for a promising type of generating profile curves comes from helical minimal surfaces. Wunderlich [18] showed that they are generated as envelope of a cylinder with a catenary as cross section and rulings orthogonal to the helical axis. Hence, we will generate helical CRPC surfaces by profile curves which are the contact curves of cylinders with rulings orthogonal to the helical axis. This fits very well to the idea of using conjugate tangents. It leads to a differential equation which can be solved explicitly, and serves as a basis for a classification of possible shapes. In particular, we analyze the singularities of positively curved helical CRPC surfaces. Their existence could only be expected, but not be shown at the discrete models by Wang and Pottmann [16].

## 2 BASIC CONCEPTS

### 2.1 Curvature Elements of CRPC Surfaces

We are interested in surfaces which possess a constant ratio $a = \kappa_1/\kappa_2$ of principal curvatures $\kappa_1, \kappa_2$. For such a CRPC surface, the Dupin indicatrices $i(p)$, given in the principal frame as

$$\kappa_1 x_1^2 + \kappa_2 x_2^2 = \pm 1, \tag{1}$$

are similar to each other. These are radial diagrams of $1/\sqrt{|\kappa_n|}$, with $\kappa_n$ as normal curvature at a tangent direction $(x_1, x_2)$. As we rule out developable surfaces, characterized by vanishing Gaussian curvature $K = 0$, we have to distinguish two cases.

(1) $K > 0$. The Dupin indicatrix $i(p)$ is an ellipse and only one sign is necessary on the right hand side of (1), depending on the sign of $\kappa_i$. The diagonals in the axis rectangle of the ellipse form the angle $\alpha = \arctan\sqrt{\kappa_1/\kappa_2}$ with the first principal direction ($x_1$-axis of the principal frame).

(2) $K < 0$. The indicatrix consists of a pair of conjugate hyperbolas. Their common asymptotes are the asymptotic directions ($\kappa_n = 0$) at angle $\alpha = \arctan\sqrt{|\kappa_1/\kappa_2|}$.

*CRPC surfaces are characterized by a constant angle $\alpha$*. They are a concept of equiform geometry, which is based on the group of Euclidean similarity transformations. Under a similarity we understand here the composition of a congruence transformation and uniform scaling.

Our approach is based on *conjugate surface tangents*. A tangent direction $(x_1, x_2)$ and its conjugate direction $(\bar{x}_1, \bar{x}_2)$, satisfy

$$\kappa_1 x_1 \bar{x}_1 + \kappa_2 x_2 \bar{x}_2 = 0 . \tag{2}$$

They describe conjugate diameters of the Dupin indicatrix. Principal directions are the only ones which are conjugate and orthogonal. Asymptotic directions are self-conjugate. The diagonal directions mentioned in the case $K > 0$ are characterized as those which are conjugate and symmetric with respect to the principal directions. They are called *characteristic directions*.

Conjugate directions have the following geometric meaning. Given a curve $c$ in a surface, we consider the envelope of surface tangent planes along $c$. This is a certain developable surface $D$. At each point of the curve $c$, the curve tangent $T$ and the ruling $\bar{T}$ of $D$ are conjugate tangents.

The bilinear conjugacy relation (2) shows that the mapping between conjugate tangents is an involutory projective map $\pi$ in the pencil of surface tangents at a surface point $p$. As maps in pencils of lines, they are congruent over the entire CRPC surface and characterized by the angle $\alpha$. For $K < 0$, the two real fixed lines of $\pi$ form the angle $2\alpha$. For $K < 0$, the pair of lines which correspond in $\pi$ and are symmetric with respect to the principal (orthogonal conjugate) directions are the characteristic directions; they also form the angle $2\alpha$.

To visualize and later apply this concept, we intersect the pencil with a circle $s$ (Steiner circle) that passes through $p$ and lies in the tangent plane $\tau(p)$ at $p$ (see Figure 1). Each tangent $T$ intersects $s$ in two points: $p$ and another point which we call $t'$. Thus the map $\pi : T \mapsto \bar{T}$ is transformed to a projective map $\pi_s : t' \mapsto \bar{t}'$ on the Steiner circle $s$. The connecting lines of corresponding points $t', \bar{t}'$ pass through a fixed point, known as involution center $I_s$. Connecting $I_s$ with the center $m_s$ of $s$ and intersecting this line with $s$ yields two points on the orthogonal conjugate (i.e. principal) tangents. For $K < 0$, $I_s$ is outside the circle $s$ and the fixed points of $\pi_s$ are the contact points of tangents from $I_s$ to $s$. These points lie in the asymptotic tangents. For $K > 0$, the line orthogonal to $I_s m_s$ intersects $s$ in 2 points which lie on the characteristic tangents. In both cases, a constant angle $\alpha$ shows that *a CRPC surface is characterized by a constant ratio between the radius $r_s$ of the Steiner circle $s$ and the distance $d_s = \|I_s - m_s\|$ of the involution center and the center of $s$* ($|\cos 2\alpha| = d_s/r_s$ for $K > 0$, $|\cos 2\alpha| = r_s/d_s$ for $K < 0$).[1]

The relation between curvature ratio $a$ and the ratio $k = d_s/r_s$ follows immediately from elementary geometry (see Figure 1, middle). Due to the equal angles over a circular arc, we can choose the center of $s$ on the first principal tangent. We find for $a > 0$ ($K > 0$),

$$k = |\cos 2\alpha| = \left|\frac{1 - \tan^2 \alpha}{1 + \tan^2 \alpha}\right| = \left|\frac{1 - a}{1 + a}\right| .$$

---

[1] $\cos 2\alpha < 0$ occurs when $I_s$, in Figure1 middle, lies on the left side of $m_s$.

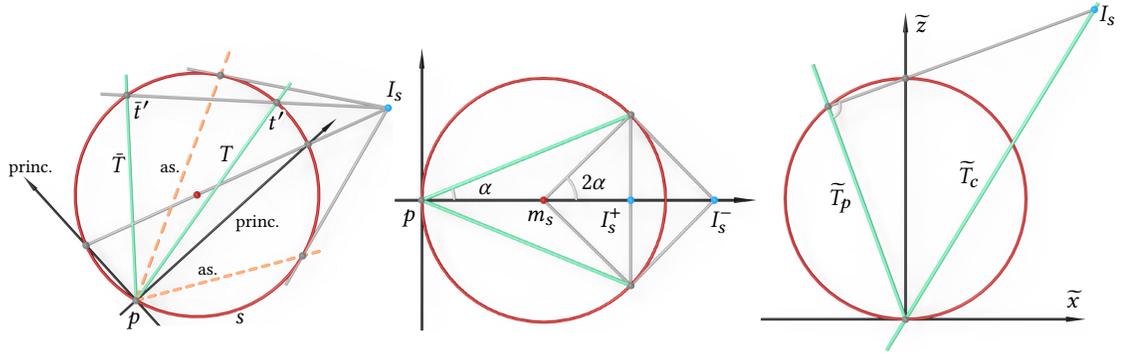

Fig. 1. Left: Involutory projective map in a pencil of lines through a point $p$, projected onto a circle $s$ through $p$: corresponding points on the Steiner circle $s$ are collinear with the involution center $I_s$. Middle: Relation between characteristic angle $\alpha$ and involution center $I_s^+$ (for $K > 0$) and $I_s^-$ (for $K < 0$). Right: Application to helical CPRC surfaces.

For $a < 0$, we have $1/k = \cos 2\alpha$ and thus in both cases obtain $k = |1-a|/|1+a|$. For most surfaces, also helical surfaces, there is no clear rule which principal direction to call the first one. Hence, it is an advantage that our parameter $k$ is agnostic to that: replacing $a$ by $1/a$, $k$ remains unchanged. Later, we will also use the parameter

$$k^2 = \left(\frac{1-a}{1+a}\right)^2 = 1 - \frac{4}{2+a+1/a}. \tag{3}$$

Let us note the special cases which we are not interested in:

- Since neither plane nor sphere are helical, it is impossible to have $k = 0, a = 1$;
- Developable surfaces, i.e. $k = 1, a \in \{0, \infty\}$;
- Minimal surfaces, i.e. $k = \infty, a = -1$.

## 2.2 Helical Surfaces.

In Euclidean $\mathbb{R}^3$, we use Cartesian coordinates $(x, y, z)$ and consider a helical motion about the $z$-axis $A$. It is composed of a uniform rotation about $A$ and a translation parallel to it. We view $A$ as vertical and call the projection onto the $xy$-plane the top view. A helical motion appears as a pure rotation in the top view.

With $(x_0, y_0, z_0)$ as coordinates in the moving system, the position in the fixed system corresponding to a rotation angle $\theta$ is

$$x(\theta) = x_0 \cos\theta - y_0 \sin\theta, \ y(\theta) = x_0 \sin\theta + y_0 \cos\theta, \ z(\theta) = z_0 + p\theta.$$

The real parameter $p$ is called the pitch. Since CRPC surfaces are a concept of equiform geometry, we could set $p = 1$. However, we want to include the rotational surface case to $p = 0$ and do not yet fix $p$.

The velocity field of the helical motion is time independent. The velocity vector $v(X)$ at a point $X = (x, y, z)$ is tangent to the helical path and equals $v(X) = (-y, x, p)$. A curve $X_0(t)$ in the moving system generates a helical surface $H : X(t, u)$. It may be generated by any of its curves as long as it is transversal to the helical paths. This leaves a lot of freedom in determining a helical surface, and one can use it to simplify the search for helical surfaces which enjoy a specific property.

In order to apply conjugacy to determine helical CRPC surfaces, we note the following well known property: At each point of a helical surface, the tangent to the helical path (velocity vector) and the direction of steepest descent in the tangent plane are conjugate to each other. This follows from the fact that the envelope of tangent planes along a helical path is a developable helical surface whose rulings are lines of steepest descent in their corresponding tangent planes.

## 3 HELICAL CRPC SURFACES

### 3.1 Setting up the Characterizing Differential Equation

Having already one simple pair of conjugate directions, we just require another one to determine the involution $\pi$ of conjugate tangents. Thus, we take as profile curve $X_0(t)$ one where the envelope $D$ of tangent planes along it is simple. We choose $D$ to be a cylinder with rulings parallel to the $x$-axis and write the positions of the tangent planes of $D$ in the fixed system in the form

$$y = tz - f(t). \tag{4}$$

Notice that we always have the tangent plane in this form since the path tangents of the helical motion are never parallel to the $xy$–plane (i.e. $t \neq \infty$).

The rulings of $D$ arise by intersection with the derivative planes $z = f'(t)$. The profile curve $X_0(t)$ along which $D$ is tangent to the helical surface $H$ is the set of those points on $D$ whose velocity vectors are tangential to $D$. This yields

$$X_0(t) = (pt, tf'(t) - f(t), f'(t)). \tag{5}$$

We call this curve the *contour* for projection parallel to the $x$-axis (or just contour for short). Its tangent,

$$X_0' = (p, tf'', f''),$$

is conjugate to the projection direction $(1, 0, 0)$. The helical path tangent vectors at the contour are

$$(f - tf', pt, p),$$

and they are conjugate to the directions of steepest descent in the tangent planes of $D$.

We rotate each tangent plane of $D$ (and $H$) about the $x$-parallel ruling of $D$ so that it becomes parallel to the $xz$-plane. The new coordinates $(\tilde{x}, 0, \tilde{z})$ of tangent vectors are related to the original ones $(x, y, z)$ by

$$(\tilde{x}, \tilde{z}) = \left(x, z\sqrt{1 + t^2}\right).$$

In the rotated position, we set up the involution of conjugate tangents. The line of steepest descent corresponds to the path tangent $\tilde{T}_p$,

$$(0, 1) \mapsto \left(f - tf', p\sqrt{1 + t^2}\right),$$

and the horizontal direction corresponds to the contour tangent $\tilde{T}_c$,

$$(1, 0) \mapsto \left(p, f''\sqrt{1 + t^2}\right).$$

As radius and center of the Steiner circle we choose $r_s = \sqrt{1 + t^2}/2$ and $m_s = (0, r_s)$ (see Figure 1, right). This yields the involution center

$$I_s = \lambda \left(p, tf''\sqrt{1 + t^2}\right),$$

with

$$\lambda = \frac{1+t^2}{t^2(1+t^2)f'' - (tf' - f)}.$$

Now we have to express a constant ratio $k \geq 0$ between its distance to $m_s$ and $r_s$,

$$\|I_s - m_s\|^2 = k^2 r_s^2,$$

which finally yields the *differential equation that characterizes helical CRPC surfaces*,

$$4p^2(1+t^2) + ((1+t^2)f'' + (tf' - f))^2 = k^2((1+t^2)f'' - (tf' - f))^2. \tag{6}$$

Due to the invariance under similarities, we can choose any pitch $p \neq 0$ for helical surfaces that are non-rotational and may then set $4p^2 = 1$. Thus in the rest of the paper, we only consider

$$(1+t^2) + ((1+t^2)f'' + (tf' - f))^2 = k^2((1+t^2)f'' - (tf' - f))^2. \tag{7}$$

**Remark 1.** *Helical minimal surfaces* are obtained with $k = \infty$ and therefore characterized by

$$(1+t^2)f'' - tf' + f = 0. \tag{8}$$

However, the derivation of this equation can be shortened a lot: The involution of conjugate tangents is the reflection at the (orthogonal) asymptotic directions. Since $x$-parallel and steepest tangent are orthogonal, the conjugate tangents, namely path tangent and contour tangent, must also be orthogonal. This yields (8), in which the pitch $p$ does no longer appear. Hence the tangent cylinders orthogonal to the axis are the same as for rotational surfaces ($p = 0$). This yields Wunderlich's generation of helical minimal surfaces as envelopes of a cylinder with a catenary as orthogonal cross section [18].

### 3.2 Solution of the Differential Equation

In this subsection we introduce an approach to find a parametric solution of ODE (7).

By setting $g = tf' - f$ we can reduce ODE (7) to

$$(1+t^2) + \left(\left(t + \frac{1}{t}\right)g' + g\right)^2 = k^2\left(\left(t + \frac{1}{t}\right)g' - g\right)^2. \tag{9}$$

Further, let $g = u\sqrt{1+t^2}$ and notice that $(t + \frac{1}{t})g' - g = \frac{u'}{t}(1+t^2)^{\frac{3}{2}}$, we get

$$(1+t^2) + \left(\frac{u'}{t}(1+t^2)^{\frac{3}{2}} + 2\sqrt{1+t^2}u\right)^2 = k^2\left(\frac{u'}{t}(1+t^2)^{\frac{3}{2}}\right)^2$$

or

$$1 + \left(2u + \left(t + \frac{1}{t}\right)u'\right)^2 = k^2\left(\left(t + \frac{1}{t}\right)u'\right)^2.$$

This leads us to a substitution for some function $s = s(t)$ such that [2]

$$\begin{cases} 2u + \left(t + \frac{1}{t}\right)u' = \sinh(s), \\ k\left(t + \frac{1}{t}\right)u' = \cosh(s). \end{cases}$$

---

[2] We may also set $k(t + \frac{1}{t})u' = -\cosh(s)$ to yield a solution on the other side of the $xz$-plane, which is equivalent to '+$\cosh(s)$' case up to a helical motion.

It immediately shows that
$$2u = \sinh(s) - \frac{\cosh(s)}{k}. \tag{10}$$

Taking the derivative
$$2u' = s'\left(\cosh(s) - \frac{\sinh(s)}{k}\right),$$

and combining with the last formula in the substitution we obtain
$$s'(k - \tanh(s)) = \frac{2t}{1+t^2}.$$

This equation is separable in variables and thus can be solved easily:
$$ks - \ln(\cosh(s)) + C' = \ln(1+t^2),$$

or
$$t^2(s) = \frac{Ce^{ks}}{\cosh(s)} - 1,$$

where $C'$ is a constant and $C = e^{C'} > 0$. Together with equation (10) we obtain two parametric solutions for ODE (9),
$$g(s) = \frac{1}{2}\left(\sinh(s) - \frac{\cosh(s)}{k}\right)\sqrt{\frac{Ce^{ks}}{\cosh(s)}}, \quad t(s) = \sqrt{\frac{Ce^{ks}}{\cosh(s)} - 1},$$

the second one is just $(g(s), -t(s))$. Finally, since $g'(t) = tf''(t)$ and by total differentiation we have
$$\frac{g'(s)}{t'(s)} = g'(t) = tf''(t) = t(s)\frac{(f'(t))'_s}{t'(s)}.$$

Thus,
$$f'(t(s)) = \int \frac{g'(s)}{t(s)} ds. \tag{11}$$

In practice, the above expressions are not easy to handle, so we prefer algebraic parametrizations
$$g(s) = \frac{(k-1)s^2 - (k+1)}{4ks}\sqrt{\frac{2Cs^{k+1}}{s^2+1}}, \quad t(s) = \sqrt{\frac{2Cs^{k+1}}{s^2+1} - 1}, \quad \left(s > 0, \ \frac{2Cs^{k+1}}{s^2+1} > 1\right). \tag{12}$$

This is just a replacement of $e^s$ by $s$. Then the expressions of the contour curves are
$$X_0(s) = \left(\frac{t(s)}{2}, g(s), \int \frac{g'(s)}{t(s)} ds\right), \tag{13}$$

$$X_1(s) = \left(-\frac{t(s)}{2}, g(s), -\int \frac{g'(s)}{t(s)} ds\right). \tag{14}$$

Here we have two constants. One arises from the integral of the $z$-coordinate, which can be ignored under a helical motion. The other one is $C$, which influences the shape of the contour and thus the shape of the helical surface.

Summarizing, we have explicitly determined all helical CRPC surfaces:

THEOREM 3.1. *Any helical surface $X$ with a constant principal curvature ratio $a$ is generated as follows. Let $k = \left|\frac{a-1}{a+1}\right|$, functions $g(s), t(s)$ be defined as in (12) and $H(v, \cdot)$ be the helical motion with pitch $1/2$,*
$$H(v, (x_0, y_0, z_0)) = (x_0 \cos v - y_0 \sin v, x_0 \sin v + y_0 \cos v, z_0 + \frac{v}{2}).$$

*Then $X$ can be parametrized as*
$$X(v, s) = 2pH(v, X_i(s)), \ s \in I_C, \ v \in \mathbb{R}, \ i = 0, 1$$

where $X_i(s)$ is one of the curves in (13) and (14), $I_C := \left\{s : s > 0, \frac{2Cs^{k+1}}{s^2+1} > 1\right\}$[3] and $p$ is the pitch of $X$.

## 3.3 Joining Two Parts to a Single $C^\infty$ Solution

Two examples of helical CRPC surfaces based on this explicit representation and using the profile curve $X_0$ are shown in Figure 2. It may seem that the two profile curves $X_0(s), X_1(s)$ generate two different CRPC surfaces. However, we will now show that this is not really the case. *The two curves generate two helical surfaces which can be joined to a single smooth (meaning $C^\infty$) CRPC surface* (Figure 3).

Of course, $X_0(s)$ and $X_1(s)$ are smooth everywhere in the domain $I_C = \left\{s : s > 0, \frac{2Cs^{k+1}}{s^2+1} > 1\right\}$. We now show that the complete contour curve is obtained by gluing curves $X_0(s)$ and $X_1(s)$ properly together (Figure 3).

To show this, we focus on the continuous function $h(s) = \frac{2Cs^{k+1}}{s^2+1}$, where $s \in [0, \infty)$. Since $h(0) = 0$ and $h(s) > 1$ for $s \in I_C$, by the Intermediate Value Theorem there exists a positive number $s_0 = \inf I_C$ such that $h(s_0) = 1$ or equivalently $2Cs_0^{k+1} = s_0^2 + 1$. No matter whether $k < 1$ or $k \geq 1$, it is easy to observe that the function $h(s)$ is strictly increasing at $s = s_0$ which gives $h'(s_0) > 0$ or equivalently $(k-1)s_0^2 + k + 1 > 0$. Since $g'(s)$ and $\sqrt{s^2+1}$ are continuous near $s_0$, when $s \to s_0^+$ we have

$$\frac{g'(s)}{t(s)} = \frac{g'(s)\sqrt{s^2+1}}{\sqrt{2Cs^{1+k} - s^2 - 1}} = O\left(\frac{1}{\sqrt{2Cs^{1+k} - s^2 - 1}}\right) = O\left(\frac{1}{\sqrt{s - s_0}}\right). \quad (15)$$

The last equality holds because

$$\lim_{s \to s_0^+} \frac{2Cs^{k+1} - s^2 - 1}{s - s_0} = \lim_{s \to s_0^+} (2C(k+1)s^k - 2s) = \frac{(k-1)s_0^2 + k + 1}{s_0} > 0.$$

Again, the last equality holds because $2Cs_0^{k+1} = s_0^2 + 1$.

The above argument guarantees the existence of the integral $\int \frac{g'(s)}{t(s)} ds$ when $s \to s_0^+$. In this regard, we may add a suitable constant to the $z$-coordinate of $X_1(s)$ s.t. $X_0(s_0) = X_1(s_0)$ and glue these two branches at this point. So the final task of this section is to show the smoothness at $s = s_0$. The idea is to treat the $y, z$-coordinates as functions of the $x$-coordinate. We start with a lemma.

LEMMA 3.2. *Consider polynomials $p(s) = \sum_{i=0}^n a_i(s)t^i(s)$ in $t(s)$ with coefficients $\{a_i(s)\}_{i=0}^n$ which are all smooth at $s = s_0$ and set*

$$P_{odd} = \{p(s) : a_{2i}(s) = 0 \ \forall i \geq 0\}, \ P_{even} = \{p(s) : a_{2i+1}(s) = 0 \ \forall i \geq 0\}.$$

*Then*

$$p(s) \in P_{even} \ (or \ P_{odd}) \implies \frac{p'(s)}{t'(s)} \in P_{odd} \ (or \ P_{even}).$$

PROOF. Since $h(s)$ is smooth at $s = s_0$ and $h'(s_0) = \frac{(k-1)s_0^2 + (k+1)}{s_0} > 0$, the function $1/h'(s)$ is smooth at $s = s_0$ as well. Assume that $p(s) \in P_{even}$(or $P_{odd}$) and notice that $t'(s) = \frac{h'(s)}{2t(s)}$. Clearly, for a monomial $a_i(s)t^i(s), i \geq 1$ of $p(s)$, we have

$$\frac{(a_i(s)t^i(s))'_s}{t'(s)} = \frac{ia_i(s)t'(s)t^{i-1}(s) + a_i'(s)t^i(s)}{t'(s)} = ia_i(s)t^{i-1}(s) + \frac{2a_i'(s)t^{i+1}(s)}{h'(s)},$$

and for $i = 0$, we have

$$\frac{(a_i(s))'_s}{t'(s)} = \frac{2a_i'(s)t(s)}{h'(s)}.$$

---
[3] We always assume that $C$ is large enough s.t. $I_C$ is nonempty.

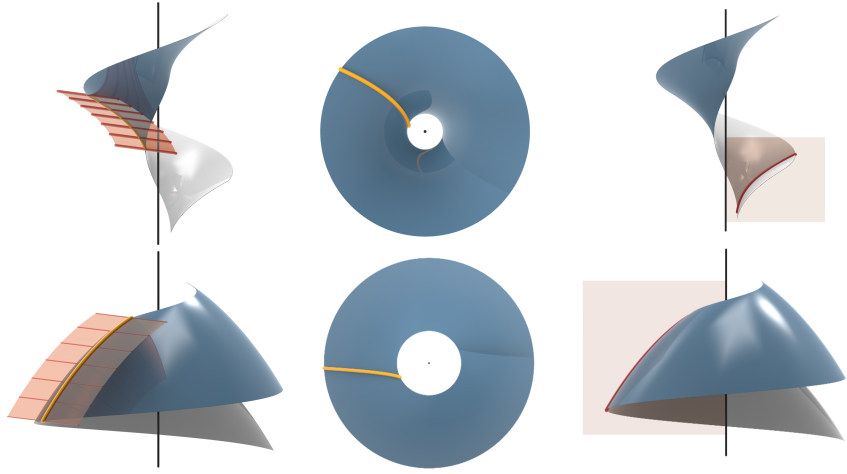

Fig. 2. Two helical CRPC surfaces generated by the profile curve $X_0$. The surface in the first row has a curvature ratio $a < 0$, the one in the 2nd row belongs to $a > 0$. 1st column: Cylinder surface tangent to the helical CRPC surface along $X_0$. 2nd column: Top view of $X_0$. 3rd column: Intersection with the $(x, z)$–plane.

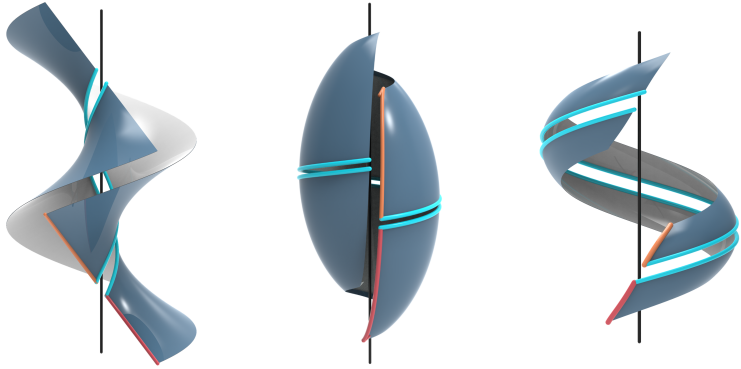

Fig. 3. The two curves $X_0$, $X_1$ given equations (13) and (14) generate the upper half and the lower half of a single smooth helical CRPC surface.

These show that the obtained monomials have degrees from $\{i-1, i+1\}$ which are different in parity from the degree of $a_i(s)t^i(s)$. Since all monomials of $p(s)$ have the same parity, we have $\frac{p'(s)}{t'(s)} \in P_{odd}$ (or $P_{even}$). □

Let us now return to the glued curve. Define two polynomial sequences $p_n(t)$ and $q_n(t)$ in $t(s)$: $p_0(t) = g(s)$ and $q_0(t) = \frac{2g'(s)}{h'(s)}$, and for $n \geq 1$

$$p_n(t) = (p_{n-1}(t))'_t = \frac{(p_{n-1}(t(s))'_s}{t'(s)}, \qquad q_n(t) = (q_{n-1}(t))'_t = \frac{(q_{n-1}(t(s))'_s}{t'(s)}.$$

Since both $p_0(t)$ and $q_0(t)$ are polynomials, free of $t(s)$, they belong to $P_{even}$. By Lemma 3.2 it immediately follows that $p_{2n}(t), q_{2n}(t) \in P_{even}$ and $p_{2n+1}(t), q_{2n+1}(t) \in P_{odd}$ for all $n \geq 0$. Thus we have

$$X_0(t) = \left(\frac{t}{2}, p_0(t), \int \frac{g'(s)}{t(s)} ds\right), \quad X_1(-t) = \left(-\frac{t}{2}, p_0(t), c - \int \frac{g'(s)}{t(s)} ds\right),$$

$$X_0^{(1)}(t) = \left(\frac{1}{2}, p_1(t), q_0(t)\right), \ X_1^{(1)}(-t) = \left(\frac{1}{2}, -p_1(t), q_0(t)\right), \tag{16}$$

$$X_0^{(n)}(t) = (0, p_n(t), q_{n-1}(t)), \ X_1^{(n)}(-t) = \left(0, (-1)^n p_n(t), (-1)^{n-1} q_{n-1}(t)\right), n \geq 2.$$

Here we should keep in mind that we treat $X_0, X_1$ as vector functions of their first coordinates (doubled), that is why we use $X_1^{(n)}(-t)$ rather than $X_1^{(n)}(t)$. The constant $c$ will make sure that $X_0(s)$ meets $X_1(s)$ at $s = s_0$. Finally notice that $t(s_0) = 0$ and $p(s_0) = 0$ whenever $p(s) \in P_{odd}$, which implies $X_0^{(n)}(0) = X_1^{(n)}(0)$ for all $n \geq 0$, i.e. the glued curve is smooth at $s = s_0$.

### 3.4 Top Views of the Profile Curves

Helical CRPC surfaces are a generalization of helical minimal surfaces. For the latter, it is known that the profiles (contours for parallel projection orthogonal to the helical axis) appear as hyperbolas in the top view [18]. Hence, it is natural to see how complicated the top views get for CRPC surfaces. We will show that the top views of the profiles of CRPC surfaces are algebraic curves for rational values of $k$.

PROPOSITION 3.3. *For any rational value of $k$, the top view of the profile curve $X_i(s)$ of a corresponding helical CRPC surface lies in an algebraic curve.*

PROOF. The top view of the profile is given by $(x, y) = (\pm t/2, g)$. Assume $k = \frac{n}{m}$ and notice that

$$t^2(s) + 1 = \frac{2Cs^{k+1}}{s^2 + 1},$$

$$g^2(s) = \left(\frac{(k-1)s^2 - (k+1)}{4k}\right)^2 \frac{2Cs^{k-1}}{s^2 + 1}.$$

Thus both $(t^2(s)+1)^m$ and $g^{2m}(s)$ are rational functions of $s$. In fact, there is a polynomial $P(t, g)$ such that $P(t(s), g(s)) = 0, \ \forall s$. More specifically, we have

$$\frac{g^2(s)}{t^2(s) + 1} = \frac{1}{16k^2}\left((k-1)^2 s^2 + \frac{(k+1)^2}{s^2} - 2(k^2 - 1)\right). \tag{17}$$

This quadratic equation in $s^2$ allows us to express $s^2$ in the form $s^2 = A + B$ where $A, \ B^2$ are rational functions of $(t, g)$. Further, let

$$F := (t^2(s) + 1)^{2m} = \frac{(2C)^{2m}(s^2)^{n+m}}{(s^2+1)^{2m}} = \frac{P_1 + Q_1 B}{P_2 + Q_2 B}$$

where $P_i, \ Q_i$ are rational functions of $(t, g)$. This step can be simplified when $m, n$ have the same parity, in which case

$$F := (t^2(s) + 1)^m = \frac{(2C)^m (s^2)^{\frac{n+m}{2}}}{(s^2+1)^m} = \frac{P_1 + Q_1 B}{P_2 + Q_2 B}.$$

Finally, since

$$B^2 = \left(\frac{P_1 - FP_2}{Q_1 - FQ_2}\right)^2$$

we get the desired polynomial by taking the numerator of

$$B^2(Q_1 - FQ_2)^2 - (P_1 - FP_2)^2 = 0.$$

□

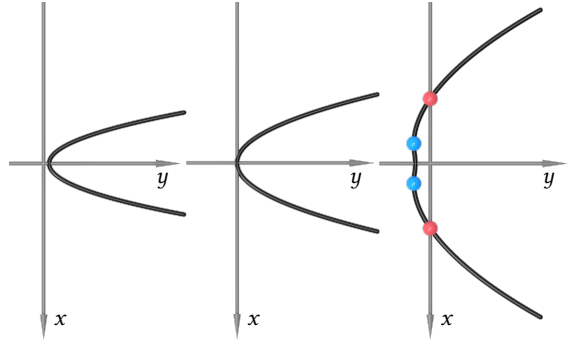

Fig. 4. Top views of the profile for $k = 3$ and different values of $C$. Left: $C = 1/8$; the contour lies on one side of the $(x, z)$–plane. Middle: $C = 3/8$; the contour touches the $(x, z)$–plane and the $z$–axis lies on the helical surface. Right: $C = 10$; the contour crosses the $(x, z)$–plane. Using the notation of Figure 1, $\tilde{T}_p$ is in asymptotic direction at the red points and $\tilde{T}_c$ is in asymptotic direction at the blue points.

By equation (17) we can identify $t^2 + 1$ as common denominator of $A$ and $B$, and thus all $P_i$, $Q_i$ admit powers of $t^2 + 1$ as their denominators. By tracking the degree in each step of the proof, it is not hard to obtain $4(3m + n)$ *as an upper bound of the degree of the final algebraic curve.*

**Example 3.1.** We compute the top view of $X_0(s)$ for $k = 3$. Based on the proof of Proposition 3.3, we have

$$\frac{g^2(s)}{t^2(s) + 1} = \frac{1}{36}\left(s^2 + \frac{4}{s^2} - 4\right).$$

So $s^2 = A + B$ where

$$A = \frac{18g^2 + 2t^2 + 2}{t^2 + 1}, \quad B^2 = \frac{36g^2(9g^2 + 2t^2 + 2)}{(t^2 + 1)^2}$$

and therefore $\frac{1}{s^2} = \frac{A-B}{4}$.

On the other hand,

$$t^2 + 1 = \frac{2Cs^4}{s^2 + 1} = \frac{2Cs^2}{1 + \frac{1}{s^2}} = \frac{2CA + 2CB}{1 + \frac{A}{4} - \frac{B}{4}}.$$

Set $P_1 = 2CA$, $Q_1 = 2C$, $P_2 = 1 + \frac{A}{4}$, $Q_2 = -\frac{1}{4}$, $F = t^2 + 1$ and take the numerator[4] of $B^2(Q_1 - FQ_2)^2 - (P_1 - FP_2)^2$ we have

$$\left(\frac{C}{3} - \frac{1}{16} - \frac{x^2}{4} - \frac{y^2}{4}\right)(4x^2 + 1)^2 - \frac{4C^2}{9}(4x^2 + 1) + 6Cy^2(4x^2 + 3y^2 + 1) = 0.$$

Here we inserted $t = 2x$, $g = y$ to obtain the top view of $X_0(s)$ and see that the curves are of algebraic order 6. They pass through the absolute points and possess contact of order 3 with the ideal line at the ideal point of the $x$-axis. Figure 4 shows the resulting curves for different values of $C$.

### 3.5 Singularities

In a study of discrete CRPC surfaces, helical ones have been computed via numerical optimization [16]. There, it appeared that the positively curved ones among them should have singularities, which however could not be computed since the formation of singularities has been prevented by fairness functionals which are required in that approach. Now, it is very easy to show that helical CRPC surfaces possess singular curves for $a > 0$.

---

[4]We recommend to use a software.

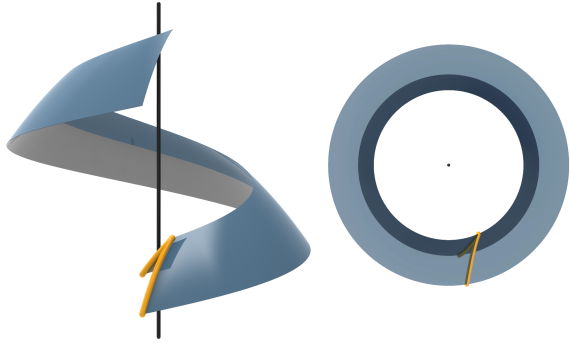

Fig. 5. Surface generated by $X_0$ near the cusp. Right: Top view.

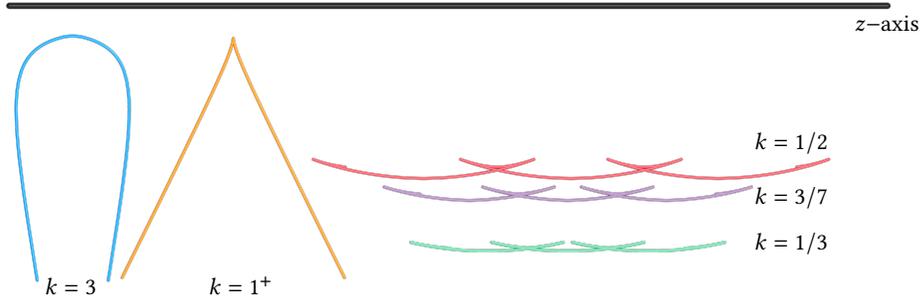

Fig. 6. Typical intersections of helical CRPC surfaces with planes through the helical axis. As $k \to 1^+$, the surface tends to a developable surface which is generated by the tangents of a helical path. As $k \to 1^-$, which is not shown in this graph, the surface tends to a cylinder, so the intersection is just a straight line parallel to the $z$−axis.

While computing $X_0'(s)$ we get a common factor of each coordinate, namely $(k-1)s^2 + 1 + k$. More precisely,

$$X_0'(s) = h'(s) \left( \frac{1}{4t(s)}, \frac{(k+1)s^2 - k + 1}{8ks\sqrt{h(s)}}, \frac{(k+1)s^2 - k + 1}{8kst(s)\sqrt{h(s)}} \right) \tag{18}$$

where

$$h'(s) = \frac{2Cs^k((k-1)s^2 + 1 + k)}{(s^2+1)^2}.$$

Thus the cusp occurs only if $s = s_k = \sqrt{\frac{1+k}{1-k}}$, in which case we must have $k < 1$ (i.e. $a > 0$); see Figure 5. In this case, $t(s)$ admits a maximum at $s = s_k$ and then decreasing till $t(s_0') = 0$ for some $s_0' > s_k > s_0$. With a similar approach we know that $X_0(s)$ and $X_1(s)$ can be glued smoothly at $s = s_0'$ as well (Figure 3). If we keep the gluing procedure with a bunch of profile curves, then we will obtain a periodic curve. For different values of $k$ the intersections of helical CRPC surfaces with the $(x, z)$−plane are shown in Figure 6.

Under the helical motion, $X_0(s_k)$ generates a singular curve which separates the helical surface into two parts. In the next subsection, we will show that is is natural to consider the ratio of principal curvatures switching to its reciprocal value when moving across a singularity. So, for convenience we split $X_0(s)$ into two curves $X_0^+(s)$, $s \geq s_k$ and $X_0^-(s)$, $s \leq s_k$.

One question arises: Do these two curves actually stop at $s = s_k$ or do they just share the same tangent here? In other words, can any of them, say $X_0^-(s)$, be continued a little bit to $s = s_k + \epsilon$ for some $\epsilon > 0$? The answer is no! Before a rigorous proof we would like to illustrate such a situation with a simple example.

**Example 3.2.** Consider the initial value problem

$$y'(x) = 1 + \sqrt{y(x)}, \ y(0) = 0.$$

By separation of variables we can solve it implicitly

$$2\sqrt{y(x)} - 2\ln\left(1 + \sqrt{y(x)}\right) = x.$$

Clearly, this function is well defined in the first quadrant because it is monotone. But it is impossible to continue this solution at the origin since $y'(0) = 1$ and any continuation will cause a negative value for $y(x)$.

The above example can be generalized to our case, i.e. ODE (7). If we set $w(t) = (t + \frac{1}{t})g'(t)$, then ODE (7) is a quadratic equation

$$(1 + t^2) + (w + g)^2 = k^2(w - g)^2$$

of $w$, which essentially consists of two equations. The two solutions meet each other when the discriminant reaches zero, i.e.

$$16k^2g^2 + 4(k^2 - 1)(1 + t^2) = 0.$$

Plugging $s = s_k$ into the above equation, we can show that $16k^2g^2(s_k) + 4(k^2 - 1)(1 + t^2(s_k)) = 0$.

On the other hand, although $X_0'(s_k) = 0$, the limit tangent at $s = s_k$ does exist. According to equation (18), the $x, y$–coordinates of a tangent vector are

$$T(s) = \left(\frac{1}{4t(s)}, \frac{(k+1)s^2 - k + 1}{8ks\sqrt{h(s)}}\right).$$

We set $D(t, g) = 16k^2g^2 + 4(k^2 - 1)(1 + t^2)$, which implies $\nabla D = (8(k^2 - 1)t, 32k^2g)$, and then we have

$$\left.\frac{\partial D}{\partial T}\right|_{s=s_k} = \nabla D(t(s_k), g(s_k)) \cdot T(s_k) = -6 - 2k^2 < 0.$$

This means that any continuation at $s = s_k$ will cause a negative value for $D(t, g)$. Thus, $X_0^-(s)$ cannot be continued at $s = s_k$. A similar argument applies to $X_0^+(s)$.

### 3.6 Classification

Finally, we aim at a classification for all helical CRPC surfaces with a given constant $k$. There is an important difference to the well-known case of rotational CRPC surfaces. There, the only parameter which influences the shape, is the constant curvature ratio $a$. In contrast, the shape of a helical CRPC surface does not only depend on $k$ (recall $k = |1 - a|/|1 + a|$), but also on the constant $C$. We will see that for negatively curved CRPC surfaces, we cannot naturally select one principal curvature direction as the first and the other as the second one. In other words, we cannot distinguish between $a$ and $1/a$. However, for positive $K$, there is a natural way to select one direction as the first and the other one as the second principal direction. In Fig. 3, these types correspond to the middle and right image. In Fig.4 the types correspond to the longer and shorter curve segments between two cusps of the intersection with a plane through the axis.

For convenience, we adopt the following basic setup. When we use the notation $X(v, s)$ from Theorem 3.1, it refers to the surface which is generated by the complete profile, obtained by gluing $X_0(s)$ and $X_1(s)$, and the pitch is always assumed to be $1/2$. It is necessary to clarify a new parametrization for $X(0, s)$. We only give a conceptual description, since the rigorous procedure is straightforward but wordy. The idea is to rewrite the curve as a function of $t = t(s)$.

(1) When $k > 1$ (i.e. $a < 0$), $t(s)$ is monotone on $[s_0, \infty)$. We set $X(0, t)$, $t \in \mathbb{R}$ s.t.
$$X(0, t) = \begin{cases} X_0(s_0) = X_1(s_0), & t = 0, \\ X_0(s^{-1}(t)), & t > 0, \\ X_1(s^{-1}(-t)), & t < 0. \end{cases}$$

(2) When $k < 1$ (i.e. $a > 0$), $t(s)$ is monotone on $[s_0, s_k]$ and $[s_k, s_0']$ respectively. Set $t_k = t(s_k)$. To avoid the cusp, we set $X^-(0, t)$, $t \in [-t_k, t_k]$ s.t.
$$X^-(0, t) = \begin{cases} X_0^-(s_0) = X_1^-(s_0), & t = 0, \\ X_0^-(s^{-1}(t)), & t \in (0, t_k], \\ X_1^-(s^{-1}(-t)), & t \in [-t_k, 0). \end{cases}$$
and $X^+(0, t)$, $t \in [-t_k, t_k]$ s.t.
$$X^+(0, t) = \begin{cases} X_0^+(s_0') = X_1^+(s_0'), & t = 0, \\ X_0^+(s^{-1}(t)), & t \in (0, t_k], \\ X_1^+(s^{-1}(-t)), & t \in [-t_k, 0). \end{cases}$$

We need to point out that in the parametrizations of $X^-(0, t)$ and $X^+(0, t)$, $X_1(s)$ is lifted to different $z$-coordinates so it can be glued with $X_0(s)$ at $s = s_0$ or $s = s_0'$, respectively. (See the last two cases of Figure 3.) Also, in $X^+(0, t)$, $X_1(s)$, $s^{-1}(t)$ stands for different inverses on the intervals $[s_0, s_k]$ and $[s_k, s_0']$ respectively. Further, since we are only interested in the surface which is free of singularities, whenever we mention a helical CRPC surface for some $k < 1$ (i.e. $a > 0$), the surface is always restricted as a part of $X^-(v, t)$ or $X^+(v, t)$. Moreover, all surfaces are translated along the $z$-axis s.t. $X(0, 0)$ (or $X^-(0, 0)$, $X^+(0, 0)$) lies on the $y$-axis. Finally, since there is the important constant $C$ in $X(v, t)$. When two or more surfaces are involved in a topic, we will use $X(v, t; C_1)$, $X(v, t; C_2)$, ... to specify the value of $C$.

*Definition 3.4.* Two helical CRPC surfaces $X(v, t; C_1)$ and $X(v, t; C_2)$ (suppose $C_1 < C_2$) are ratio-equivalent (*R-equivalent*) if there exists a one-parameter transformation $T(v, t; C)$, $C \in [C_1, C_2]$ s.t. $T(v, t; C)$ is a helical CRPC surface for all $C \in [C_1, C_2]$ and $T(v, t; C_1) = X(v, t; C_1)$, $T(v, t; C_2) = X(v, t; C_2)$.

In other words, $X(v, t; C_1)$ is continuously deformed to $X(v, t; C_2)$ by $T$ and this transformation preserves the helical CRPC property. Obviously, $T(v, t; C)$ could be $X(v, t; C)|_{[C_1, C_2]}$.

*R-equivalence for negative curvature.* The expressions for $X(v, t; C)$ immediately show that for $k > 1$, all $X(v, t; C)$'s are R-equivalent. In contrast to the rotational case, it is a little bit surprising that *one cannot distinguish between helical CRPC surfaces with negative ratios $a$ and $1/a$ (i.e. $k > 1$), respectively*. This does not mean that they all look the same. Figure 7 shows essentially different shapes, to be discussed later in this section.

*The case of positive curvature.* For $k < 1$, we have to split at singularities and here all $X^-(v, t; C)$'s or all $X^+(v, t; C)$'s are R-equivalent. However, $X^-(v, t; C)$ is not R-equivalent to $X^+(v, t; C)$. Otherwise there would be $C_1 < C_2$ s.t. $X^-(v, t; C_1)$ is R-equivalent to $X^+(v, t; C_2)$ under a transformation $T(v, t; C)$. $T$ would automatically induce a 'sub-transformation' between profile curves $X^-(0, t; C_1)$ and $X^+(0, t; C_2)$. This is impossible, since we already mentioned in the last subsection that $X^-(0, t; C_1)$ and $X^+(0, t; C_2)$ are constructed by two different ODEs which do not share any common solutions.

*Distinguishing principal directions.* The principal directions of a rotational surface are the tangents of its parallel circles and meridians. Thus, the associated principal curvatures can be labeled as $\kappa_1$ and $\kappa_2$ respectively. If we set $\kappa_1/\kappa_2 = a$ or $\kappa_1/\kappa_2 = 1/a$ for some constant $a \neq 0, \pm 1$, there will be two essentially different (i.e. non-similar) surfaces corresponding to each case (see e.g. [16]). However, up to now everything we discussed about helical CRPC surfaces is determined by the constant $k$, where the cases $\kappa_1/\kappa_2 = a$ and $\kappa_1/\kappa_2 = 1/a$ are not distinguished at all. We will now see that *for $k < 1$, i.e. $a > 0$, there is a natural way to label the principal curvatures of a helical CRPC surface so that we can determine their ratio as either $a$ or $1/a$.*

The idea is to consider the tangents of $X(v, t)$ at the point $X(0, 0)$ where the two solutions $X_0, X_1$ have been glued together. This point generates a helical path along which the CRPC surface is tangent to a co-axial rotational cylinder (see Fig. 3). The profile curve $X(0, t)$ and the helical path $X(v, 0)$ are symmetric with respect to the $y$-axis, on which $X(0, 0)$ lies. We will now see that this point is crucial to the classification. There, path tangent $X_v(0, 0)$ and profile tangent $X_t(0, 0)$ help us to find a desired order of the principal curvatures.

We get started by another illustration of the in-equivalence of $X^-(v, t; C_1)$ and $X^+(v, t; C_2)$. Following equations (16) and (18,) we have

$$X_v^-(0, 0) = \left(-\frac{(k-1)s_0^2 - (k+1)}{4ks_0}, 0, \frac{1}{2}\right),$$

$$X_t^-(0, 0) = \left(\frac{1}{2}, 0, \frac{(k+1)s_0^2 - k + 1}{4ks_0}\right),$$

$$X_{vv}^-(0, 0) = \left(0, -\frac{(k-1)s_0^2 - (k+1)}{4ks_0}, 0\right),$$

$$X_{tt}^-(0, 0) = \left(0, \frac{(k+1)s_0^2 - k + 1}{4ks_0}, 0\right).$$

Clearly, the surface normal at $X^-(0, 0)$ is $n = (0, 1, 0)$. So, coefficients $E, G$ of the first fundamental form and $L, N$ of the second fundamental form are

$$E = X_v^- \cdot X_v^- = \frac{(s_0^2 + 1)((k-1)^2 s_0^2 + (k+1)^2)}{16k^2 2s_0^2},$$

$$G = X_t^- \cdot X_t^- = \frac{(s_0^2 + 1)((k+1)^2 s_0^2 + (k-1)^2)}{16k^2 2s_0^2},$$

$$L = n \cdot X_{vv}^- = -\frac{(k-1)s_0^2 - (k+1)}{4ks_0},$$

$$N = n \cdot X_{tt}^- = \frac{(k+1)s_0^2 - k + 1}{4ks_0}.$$

Let $\kappa_v^-$, $\kappa_t^-$ be the normal curvatures of $X^-(v, 0)$, $X^-(0, t)$ respectively. Then at $X^-(0, 0)$ we have

$$\frac{\kappa_v^-}{\kappa_t^-} = \frac{L/E}{N/G} = \frac{LG}{NE} = \frac{(1+k)s_0^4 + O(s_0^2) + (1-k)}{(1-k)s_0^4 + O(s_0^2) + (1+k)}.$$

If we let $C \to +\infty$, $X^-(0,0)$ will tend to $-\infty$ along the $y$-axis. Since we have fixed the pitch at $1/2$, moving with $X^-(0,0)$ to infinity, the surface becomes locally like a rotational one. The parameter curves $X^-(v,0)$, $X^-(0,t)$ are locally like parallel circle and meridian. So the ratio of principal curvatures can be approximated by $\kappa_v^-/\kappa_t^-$. Similarly, the ratio for $X^+(v,t)$ can be approximated by

$$\frac{\kappa_v^+}{\kappa_t^+} = \frac{(1+k)s_0'^4 + O(s_0'^2) + (1-k)}{(1-k)s_0'^4 + O(s_0'^2) + (1+k)}.$$

Consider the way how we choose $s_0$ and $s_0'$: The function $h(s) = \frac{2Cs_0^{1+k}}{s_0^2+1}$, $k < 1$ is firstly increasing and then decreasing on $[s_0, s_0']$ where $h(s_0) = h(s_0') = 1$. It is not hard to see that $C \to +\infty$ implies $s_0 \to 0$ and $s_0' \to +\infty$. Thus

$$\frac{\kappa_v^-}{\kappa_t^-} \to \frac{1-k}{1+k}, \quad \frac{\kappa_v^+}{\kappa_t^+} \to \frac{1+k}{1-k},$$

showing again that $X^-(v,t;C_1)$ and $X^+(v,t;C_2)$ are not R-equivalent.

Let us review the above discussion. It actually gives a natural way to label the principal curvatures (including the case $a < 0$): For any surface $X(v,t;C)$, suppose we have an unlabeled principal frame at $X(0,0;C)$. By increasing the value of $C$, this point is pushed to infinity and one of its principal directions continuously drives to the direction of $X_v(0,0;C)$. We label the principal curvature in this direction by $\kappa_1$ and the other by $\kappa_2$.

In the case $a < 0$, we can just copy all the calculations from $X^-(v,t)$ and obtain

$$\frac{\kappa_v}{\kappa_t} = \frac{(1+k)s_0^4 + O(s_0^2) + (1-k)}{(1-k)s_0^4 + O(s_0^2) + (1+k)}.$$

However, it will not help us to classify the surfaces. This is because $h(s)$ is increasing on $[s_0, +\infty)$ where $h(s_0) = 1$ and this interval always exists as long as $C > 0$. When $C$ starts from $0^+$ and tends to $+\infty$, one can show that $s_0$ is decreasing from $+\infty$ to $0^+$ and $X(0,0;C)$ is moving from $+\infty$ to $-\infty$ along the $y$-axis. Thus

$$\lim_{C \to +\infty} \frac{\kappa_v}{\kappa_t} = \frac{1-k}{1+k}, \quad \lim_{C \to 0^+} \frac{\kappa_v}{\kappa_t} = \frac{1+k}{1-k}.$$

In other words, the principal frame at $X(0,0;C)$ will rotate by 90 degrees when it travels through the whole $y$-axis.

*Shapes for negative curvature.* In order to classify the case $a < 0$, we have to try another approach. One way is to consider the intersection of the surface and the $(y,z)$-plane, which we call $(y,z)$-profile. This profile generates the same surface as $X(0,t)$ does. But for different values of $C$, the $(y,z)$-profiles can be classified into three major cases (see Figures 4 and 7). Recall the middle picture of Figure 4, which is related to the case $s_0 = \sqrt{\frac{k+1}{k-1}}$. Or equivalently, since $h(s_0) = 1$,

$$C = C_k := \frac{k(k-1)^{\frac{k-1}{2}}}{(k+1)^{\frac{k+1}{2}}}.$$

When $C \in (0, C_k)$ i.e. $s_0 > \sqrt{\frac{k+1}{k-1}}$, $g(s)$ is always positive. As we can see from Figure 4 left, the top view of $X_0(s)$ lies in the first quadrant. To get the upper half of the $(y,z)$-profile, we need to move each point $X_0(s) = \left(\frac{t(s)}{2}, g(s), \int \frac{g'(s)}{t(s)} ds\right)$ along its helical path until it reaches the $(y,z)$-plane. The corresponding coordinate on the $(y,z)$-profile is

$$\left(0, \sqrt{\frac{t(s)^2}{4} + g(s)^2}, \int \frac{g'(s)}{t(s)} ds + \frac{1}{2} \tan^{-1} \frac{t(s)}{2g(s)}\right).$$

Similarly, the point $\left(-\frac{t(s)}{2}, g(s), -\int \frac{g'(s)}{t(s)} ds\right)$ on $X_1(s)$ must be moved backwards to the $(y, z)$–plane, and the corresponding coordinate on the lower half of the $(y, z)$–profile is

$$\left(0, \sqrt{\frac{t(s)^2}{4} + g(s)^2}, -\int \frac{g'(s)}{t(s)} ds - \frac{1}{2} \tan^{-1} \frac{t(s)}{2g(s)}\right).$$

If we fix $X(0, 0)$ on the $y$–axis, then this $(y, z)$–profile is smooth and symmetrical with respect to the $y$–axis (Figure 7 left).

When $C > C_k$ i.e. $s_0 < \sqrt{\frac{k+1}{k-1}}$ (Figure 4 right), $X(0, 0)$ lies on the negative part of the $y$–axis. Now the points on $X_0(s)$ should be moved backwards and the points on $X_1(s)$ should be moved forward. $X_0(s)$ contributes half of the $(y, z)$–profile. Since the curve crosses the $(x, z)$–plane at $s = \sqrt{\frac{k+1}{k-1}}$ (i.e. $g(s) = 0$), the expressions are given piecewise:

$$\left(0, -\sqrt{\frac{t(s)^2}{4} + g(s)^2}, \int \frac{g'(s)}{t(s)} ds + \frac{1}{2} \tan^{-1} \frac{t(s)}{2g(s)}\right), \quad s_0 < s < \sqrt{\frac{k+1}{k-1}},$$

$$\left(0, -\sqrt{\frac{t(s)^2}{4} + g(s)^2}, \int \frac{g'(s)}{t(s)} ds - \frac{\pi}{2} - \frac{1}{2} \tan^{-1} \frac{2g(s)}{t(s)}\right), \quad s \geq \sqrt{\frac{k+1}{k-1}}.$$

The complete $(y, z)$–profile is obtained by reflecting the above curve at the $y$–axis. In this case, the $(y, z)$–profile possesses a self-intersection at the $y$–axis (Figure 7 right). This is shown as follows. The highest degree of $s$ in $g'(s)$ and $t(s)$ is the same. This means $\frac{g'(s)}{t(s)}$ tends to a (positive) constant when $s \to +\infty$. Since $\tan^{-1} \frac{t(s)}{2g(s)}$ is bounded, $\left(\int \frac{g'(s)}{t(s)} ds - \frac{1}{2} \tan^{-1} \frac{t(s)}{2g(s)}\right) \to +\infty$ when $s \to +\infty$. On the other hand, one can show that the derivative $\left(\int \frac{g'(s)}{t(s)} ds - \frac{1}{2} \tan^{-1} \frac{t(s)}{2g(s)}\right)'_s \to -\infty$ when $s \to s_0^+$. These calculations guarantee that the curve will cross the $y$–axis, i.e. possess a self-intersection.

For the last case $C = C_k$ i.e. $s_0 = \sqrt{\frac{k+1}{k-1}}$ (Figure 4 middle), $X(0, 0)$ coincides with the origin and the tangent at this point lies on the $(x, z)$–plane. This implies $\lim_{s \to s_0^+} \tan^{-1} \frac{t(s)}{2g(s)} = \frac{\pi}{2}$. So for convenience, we consider the $(x, z)$–profile instead, which is tangent to $X(0, t)$ at the origin. Again, $X_0(s)$ contributes half part of it, which is

$$\left(\sqrt{\frac{t(s)^2}{4} + g(s)^2}, 0, \int \frac{g'(s)}{t(s)} ds - \frac{1}{2} \tan^{-1} \frac{2g(s)}{t(s)}\right).$$

The other half from $X_1(s)$ is

$$\left(-\sqrt{\frac{t(s)^2}{4} + g(s)^2}, 0, -\int \frac{g'(s)}{t(s)} ds + \frac{1}{2} \tan^{-1} \frac{2g(s)}{t(s)}\right).$$

Together they form an odd (smooth) function in the $(x, z)$–plane (Figure 7 middle). Clearly, the $(y, z)$–profile is congruent to it and obtained by applying the helical motion with an angle of 90 degrees.

### 3.7 Future Work

One direction for future research is the determination of all spiral CRPC surfaces, as extensions of the known spiral minimal surfaces [19]. This is a natural question, since CRPC surfaces are invariant under Euclidean similarities and spiral surfaces are generated by one parameter subgroups of the Euclidean similarity group. The present approach based on profiles for projection orthogonal to the spiral axis should be promising, since these profiles also appear in Wunderlich's spiral minimal surfaces [19]. However, the characterizing ODE is much more complicated and thus we leave its solution for future work.

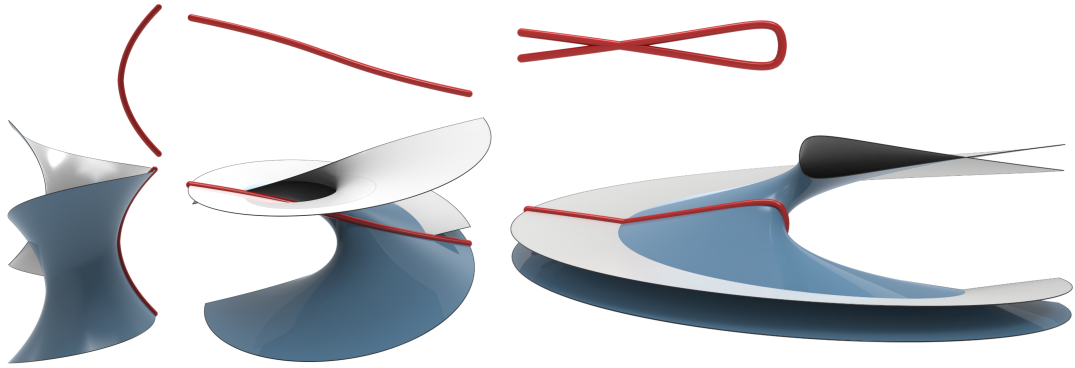

Fig. 7. The shapes of CRPC surfaces to a curvature ratio $a < 0$ (here $k = 3$) depend on the value of $C$. Left: $C = 0.01$. Center: $C = C_k = 3/8$. Right: $C = 1$.

The determination of explicit representations for CRPC surfaces apart from the so far mentioned ones is a bigger challenge, as it amounts to the solution of a rather complicated nonlinear partial differential equation. A geometric construction of general CRPC surfaces can be based on a Christoffel-type transformation of the Gaussian spherical image [2], but it remains open whether this is a path towards explicit parameterizations of CRPC surfaces.